\documentclass[12pt]{article}
\usepackage[T1]{fontenc}
\usepackage[english]{babel}
\usepackage{mathrsfs}
\usepackage{amssymb}	
\usepackage{amsthm}
\usepackage{xcolor}
\usepackage{amsmath}
\usepackage{multirow}
\usepackage{nccmath}
\usepackage{rotating}
\usepackage{subcaption}
\usepackage{float} 
\usepackage{hyperref}
\usepackage{multicol}

\usepackage{tabularx}

\newcolumntype{Y}{>{\centering\arraybackslash}X}
\usepackage{array}
\usepackage{tikz}
\usetikzlibrary{arrows.meta}
		
\usepackage[left = 2cm, right = 2cm, bottom = 2cm, top = 2cm]{geometry}

\DeclareMathOperator{\tr}{tr}
\newcommand{\dd}{{\operatorname{d}}}

\newcommand{\R}{\mathbb R}
\newcommand{\Id}{\operatorname{Id}}

\newcommand{\weg}[1]{}

\theoremstyle{plain}
\newtheorem{Theorem}{Theorem}

\newtheorem{Definition}{Definition} 
\newtheorem{Remark}{Remark}
\newtheorem{Proposition}{Proposition}

\usepackage{graphicx}
\graphicspath{ {./images/} }

\title{Singularities of Separation Coordinates}
\author{Dinmukhammed Akpan\footnote{Institute of Mathematics, Friedrich Schiller University Jena, 07743 Jena, Germany, Institute of Mathematics and Mathematical Modeling, Almaty, Kazakhstan, dinmukhammed.akpan@uni-jena.de}}

\date{}

\begin{document}

\maketitle
\begin{abstract}
We study singularities of separation coordinates
for geodesic flows of real-analytic two-dimensional
(pseudo-)Riemannian metrics. We obtain a complete local
classification of these singularities using the Nijenhuis
operator associated with a quadratic first integral.
The singular set is locally a smooth curve, an isolated
point, or two transversely intersecting curves.
We also show that orthogonal separation near one point
of an analytic disk extends to a quadratic first integral
on the entire disk and provide explicit separation
coordinates near singularities.
\end{abstract}


\tableofcontents

\section{Introduction}

We study the singularities of separation coordinates for
geodesic flows of real-analytic two-dimen\-sional metrics.  Separation of variables
is one of the central mechanisms by which the equations of classical
and quantum mechanics become explicitly tractable: it transforms the
Hamilton--Jacobi equation, and, after quantization, the corresponding
Schr\"odinger equation, into a family of lower-dimensional equations.
Its origins lie in the work of Hamilton and Jacobi; its geometric
theory was developed through the contributions of St\"ackel,
Levi-Civita, Eisenhart, and their successors; see
\cite{Benenti1991,Benenti2016}.  The question addressed in this paper is which singularities separation coordinates may have. We give a complete answer in the
two-dimensional real-analytic category.

Defining a singular point of separation coordinates requires some
care. A particular choice of coordinate functions may fail to extend
regularly even when an equivalent system of separation coordinates
is regular. We therefore
allow independent reparametrisations of the separation coordinates
before testing whether they extend across the point. A point is
singular if no such reparametrisation gives a regular extension;
see Definition~\ref{def:singular-separation-point}.
Proposition~\ref{prop:singular-point-equivalence} characterises these
points invariantly by the coincidence of the eigenvalues of the
operator associated with the separation.

For Riemannian metrics, the singularities of separation
coordinates are already understood: up to a local diffeomorphism,
there are two types, the centre of polar coordinates and a focus of
elliptic coordinates~\cite{bmf}; see Figures~\ref{fig:fourth}
and~\ref{fig:third}, respectively. Both coordinate systems are widely
used in mathematical physics. Polar coordinates separate the planar
Kepler--Coulomb problem, while elliptic coordinates are adapted to
the planar two-centre Coulomb problem. Their singularities reflect
the geometry near the centres of these systems.

The pseudo-Riemannian (Lorentzian) situation is
substantially richer.  The relevant operator may pass between real
and complex-conjugate eigenvalue regimes, and the singular set can
contain curves as well as isolated points.  Such phenomena occur
naturally in relativity: one of our basic models is the Rindler
coordinate system, which is adapted to uniformly accelerated
observers in Minkowski spacetime
\cite{wiki:rindler}; see
Figure~\ref{fig:rindler}.  In this setting, the
singularities are not exceptional artefacts of a particular
coordinate choice: they reflect genuine changes in the algebraic type
of the underlying geometric structure.

Theorem~\ref{thm:Killing-separation} describes the consequences
of nonorthogonal separation: the metric
admits a nontrivial Killing vector field, and  orthogonal
separation exists near almost every point unless this field
is identically null. In the identically null case, the metric
and the Killing field admit an explicit flat local form.
Theorem~\ref{thm:stackel} gives a criterion for orthogonal separation.
Theorem~\ref{main1} then classifies the associated operators
near scalar singular points. Together with the
$\operatorname{gl}$-regular classification in~\cite{BKM-3},
it provides all local operator forms needed below.

Theorem~\ref{thm:singular-set} shows that orthogonal separation
near one point of a real-analytic disk yields a quadratic
first integral on the entire disk. Orthogonal separation
then exists near almost every point, and its singular set
has only three possible local topological types:
a smooth curve, an isolated point, or two transversely
intersecting curves.
Finally, Theorem~\ref{thm:separation-coordinates} gives
explicit local separation coordinates near the singular
points.

\textbf{Acknowledgment.} 
The author is grateful to V.~Matveev, A.~Bolsinov and M.~Manenti for their many helpful suggestions. This research was supported by the Ministry of Science and Higher Education of the Republic of Kazakhstan (grant No. AP23483476) and by the DFG, project No.~529233771.

\section{Preliminaries and Main results}

Let $(M^2, g)$ be a (pseudo-)Riemannian surface.  On the cotangent
bundle $T^*M^2$, endowed with its canonical symplectic structure,
consider the  Hamiltonian
\[
H = \frac{1}{2}g^{-1}(p,p).
\]
In local canonical coordinates $(x,y,p_x,p_y)$, the resulting
Hamiltonian system is the geodesic flow of \(g\).  The Hamilton--Jacobi equation for the geodesic flow is
\[
H\left(
x,y,
\frac{\partial W}{\partial x},
\frac{\partial W}{\partial y}
\right)
=
c_1.
\]
A solution $W=W(x,y;c_1,c_2)$, defined on $U\times A$,
where $U\subset M^2$ and $A\subset\R^2$,
is called a \emph{Hamilton–Jacobi function} or \emph{complete integral} if it satisfies
the nondegeneracy condition
\begin{equation}
\label{eq:complete-integral-nondegeneracy}
\det\left(
\frac{\partial^2W}{\partial x^i\partial c_j}
\right)_{i,j=1}^2\ne0,
\qquad x^1=x,\quad x^2=y.
\end{equation}

The geodesic Hamilton--Jacobi equation always admits local
complete integrals near nonzero momenta. 

\begin{Definition}
The geodesic flow is said to admit
\emph{additive separation of variables} on
\(U\subset M^2\) if there exist local coordinates
\((x,y)\) on \(U\) and a complete integral of the form
\begin{equation*}
W(x,y;c_1,c_2)
=
W_1(x;c_1,c_2)
+
W_2(y;c_1,c_2)
\end{equation*}
satisfying~\eqref{eq:complete-integral-nondegeneracy}.
The coordinates $(x,y)$ are called {real separation coordinates}. 
\end{Definition}

If \((x,y)\) are separation coordinates, then independent local reparametrizations
\[
\widetilde x=\varphi(x),
\qquad
\widetilde y=\psi(y),
\qquad
\varphi'(x)\psi'(y)\neq0,
\]
as well as the interchange of \(x\) and \(y\), preserve the separated
form. We shall
refer to these transformations as \emph{admissible transformations
of separation coordinates}.

\begin{Definition}
The geodesic flow is said to admit
\emph{complex additive separation of variables} on
\(U \subset M^2\) if it possesses a real-valued complete integral of the form
\begin{equation*}
W(z,\overline z;c_1,c_2)
=
W_1(z;c_1,c_2)
+
W_2(\overline z;c_1,c_2),
\end{equation*}
where \(z\colon U\to\mathbb C\) satisfies
$\dd\operatorname{Re}z\wedge\dd\operatorname{Im}z\ne0$,
and $c_1,c_2$ are real parameters. 
The pair
\((z,\overline z)\) is called a system of
{complex-conjugate separation coordinates}.
\end{Definition}

Since the complete integral $W$ is real-valued, after
adding functions of the parameters to the separated summands,
they can be chosen complex-conjugate:
\[
W_2(\overline z;c_1,c_2)
=
\overline{W_1(z;c_1,c_2)},
\qquad
W=2\operatorname{Re}W_1(z;c_1,c_2).
\]

The complex separated form is preserved under conjugate holomorphic
reparametrizations
\[
\widetilde z=\varphi(z),
\qquad
\overline{\widetilde z}
=
\overline{\varphi(z)},
\qquad
\varphi'(z)\neq0,
\]
as well as by the interchange of \(z\) and \(\overline z\).
Again, we shall
refer to these transformations as \emph{admissible transformations
of separation coordinates}.

\subsection{Orthogonal and nonorthogonal separation of variables}

Let $H=\frac12 g^{-1}(p,p)$ be the Hamiltonian
of the geodesic flow of the metric $g$.  The geodesic flow is said to be
\emph{quadratically integrable} if there exists a function
\[
F = a(x,y)p_x^2 + b(x,y)p_x p_y + c(x,y)p_y^2
\]
such that $H$ and $F$ are functionally independent and Poisson
commute, $\{H,F\}=0$.

Separation of variables for two-dimensional geodesic flows is
closely related to the existence of first integrals polynomial
in the momenta. Nonorthogonal separation
implies the existence of a nontrivial linear in momenta first integral,
or equivalently, a nontrivial Killing vector field.
Orthogonal separation yields
a real quadratic first integral functionally independent of
the energy.

\begin{Definition}
A system of real separation coordinates $(x,y)$ is called
\emph{orthogonal} if $g(\partial_x,\partial_y)=0$.
A system of complex-conjugate separation coordinates $(z,\bar z)$
is called \emph{orthogonal} if
$g(\partial_z,\partial_{\bar z})=0$, using the complex-bilinear
extension of $g$. Otherwise, the separation is called
\emph{nonorthogonal}.
\end{Definition}

\begin{Theorem}
\label{thm:Killing-separation}
Let $g$ be a real-analytic metric near a point $P\in M^2$.
Suppose that its geodesic flow admits real or complex nonorthogonal separation of variables near $P$.
Then, in a sufficiently small connected neighbourhood $U(P)$,
there exists a nontrivial real-analytic Killing vector field $K$,
and the following assertions hold.
\begin{enumerate}
\item
If $g(K,K)\not\equiv0$, then the geodesic flow admits
real orthogonal separation of variables almost everywhere on $U(P)$.

\item
If $g(K,K)\equiv0$, then $K$ has no zeros and near every point of $U(P)$, there exist real-analytic
coordinates $(x,y)$ in which
\[
g=2\,\dd x\,\dd y,\qquad K=\partial_y.
\]
\end{enumerate}
\end{Theorem}

We now turn to orthogonal separation.
Suppose that the geodesic flow of a two-dimensional metric \(g = g_{ij}\dd x^i\dd x^j\)
admits  a quadratic first integral $F=F^{ij}(x)p_ip_j.$ 
Define the tensor field of type \((1,1)\) by

\begin{equation} \label{operator}
        L=\operatorname{adj} (g) \cdot \operatorname{adj}(F),
\qquad L^i{}_j=(g_{ab}F^{ab})\delta^i_j-F^{ik}g_{kj},
\end{equation}
where $\operatorname{adj}$ denotes the adjugate\footnote{$A \cdot \operatorname{adj}(A) = \operatorname{adj}(A) \cdot A = \det(A) \Id$} matrix and
$\delta^i_j$ is the Kronecker delta. We call \(L\) the \emph{operator associated with the pair \((g,F)\)}. Note that $L$ is self-adjoint with respect to \(g\).   The operator \(L\) is Nijenhuis; see
\cite{BolsinovMatveev2003,Akpan1} and
\cite[Section~1.1]{nijapp5}.

\begin{Theorem}\label{thm:stackel}
Let $g$ be a real-analytic metric in some neighbourhood of a point $P\in M^2$.
The geodesic flow of $g$ admits real or
complex orthogonal separation of variables near $P$
if and only if there exists a real quadratic first integral
$F$, defined near $P$, such that the associated operator
$L$, given by~\eqref{operator}, has two distinct eigenvalues
at $P$.
\end{Theorem}

\begin{Remark}
Genuinely complex orthogonal separation associated with a nonreal
complex-conjugate pair of eigenvalues of \(L\) can occur only for a
pseudo-Riemannian (Lorentzian) metric.
\end{Remark}

As an example of complex orthogonal separation, consider
\[
g=\frac{\overline{h(z)}-h(z)}{8}
\bigl(\dd z^2-\dd\overline z^2\bigr),
\qquad
F=
4\,\frac{
\overline{h(z)}\,p_z^2-h(z)\,p_{\overline z}^2
}{
\overline{h(z)}-h(z)
},
\]
where $h$ is holomorphic, $\operatorname{Im}h(z)\ne0$ and
the associated operator is $L = \operatorname{diag} \left( \frac{h(z)}{2}, \frac{\overline{h(z)}}{2} \right)$
in the coordinates $(z,\overline z)$.
Its eigenvalues are distinct and nonreal, so
Theorem~\ref{thm:stackel} gives complex
orthogonal separation.

\subsection{The associated operator near a scalar point}

As we have seen, the operator \(L\) plays a central role in the theory
of separation of variables.  We now fix
the terminology that will be used below.  Since all questions under
consideration are local, the relevant properties of \(L\) must always
be understood in a neighbourhood of the point under study.

\begin{Definition}
Let $L$ be two-dimensional an operator field near a point $P$ of a surface.
The point $P$ is called \emph{scalar} if
$L(P)=\lambda_0\Id$ for some $\lambda_0\in\R$, and
\emph{$\operatorname{gl}$-regular} otherwise.
The operator $L$ is called \emph{generic} if the algebraic type of $L$ is
constant near $P$, and \emph{singular} otherwise.
\end{Definition}

In dimension two, we distinguish the following cases.

\textbf{\(\operatorname{gl}\)-regular (both generic and singular).}
It is known that every \(\operatorname{gl}\)-regular Nijenhuis operator
is associated with a quadratically integrable metric \cite{nijapp5}.  The
classification of two-dimensional \(\operatorname{gl}\)-regular
Nijenhuis operators is given in \cite{BKM-3}.

\textbf{Generic scalar.}
In this case $L=\lambda(x,y)\Id$, so~\eqref{operator} gives
$F=2\lambda(x,y)H$.
The identity $\{H,F\}=0$ implies that $\lambda$ is constant.

\textbf{Singular scalar.}
Here $L(P)=\lambda_0\Id$, but $L$ is not scalar throughout
any neighbourhood of $P$.
Theorem~\ref{main1} gives a complete description of this case.

We impose no restriction on the signature of the metric: in dimension
two, \(g\) may be Riemannian or Lorentzian.  The displayed metrics
realizing SS5 and SS6 below are Riemannian, whereas those realizing
SS1--SS4 are Lorentzian.

\begin{Theorem} \label{main1}
Let \(L\) be a two-dimensional real-analytic operator field in a
neighbourhood of a singular scalar point \(P\in M^2\), and suppose
that \(L(P)=\lambda_0\operatorname{Id}\), where
\(\lambda_0\in\mathbb R\).  Then \(L\) is associated with a
two-dimensional metric and a nontrivial quadratic first integral of
its geodesic flow if and only if, after a constant
normalisation of the integral, there exist local coordinates centred
at $P$ in which $L$ takes
one of the following forms:

        \begin{itemize}
            \item[(SS1)] $L = \lambda_0 \operatorname{Id} + \begin{bmatrix}
                2x & -y \\
                y & 0
            \end{bmatrix}$;
            \item[(SS2)] $L = \lambda_0 \operatorname{Id} + \begin{bmatrix}
        -y(1 + (s-3)xy^{s-2}) & -x(2 + (s-1)(s-3)xy^{s-2}) \\
        y^s & -y(1 - (s-1)xy^{s-2})
    \end{bmatrix}$, $s\in\mathbb Z$, $s\geq2$;
            \item[(SS3)] $L = \lambda_0 \operatorname{Id} + \begin{bmatrix}
       x y & x^2 \\
       y^2 & xy
    \end{bmatrix}$;
            \item[(SS4)] $L = \lambda_0 \operatorname{Id} + \begin{bmatrix}
                \frac{x}{2} & 0 \\
                y & \frac{x}{2}
            \end{bmatrix}$;
            \item[(SS5)] $L = \lambda_0 \operatorname{Id} + \begin{bmatrix}
                2x & y \\
                y & 0
            \end{bmatrix}$;
            \item[(SS6)] $L = \lambda_0 \operatorname{Id} + \begin{bmatrix}
       x^2 & xy \\
       xy & y^2
    \end{bmatrix}$.
        \end{itemize}
\end{Theorem}

These normal forms describe the operators at scalar singular
points. Away from the set where the eigenvalues coincide, their
eigendirections determine the corresponding systems of orthogonal
separation coordinates, up to admissible
transformations (see Theorem \ref{thm:separation-coordinates}).

\subsection{Separation coordinates and singularities}

\begin{Definition}\label{def:singular-separation-point}
Let \((x,y)\) be separation coordinates defined on an open
set \(\Omega\subset M^2\), and let \(P\in\overline{\Omega}\).
The point \(P\) is called a {regular point of the separation
coordinates} if there exist a neighborhood \(V\) of \(P\) and an
admissible transformation of the separation coordinates, defined on
\(\Omega\cap V\), such that the transformed coordinates extend to a regular system of separation coordinates on \(V\).
Otherwise, \(P\) is called a {singular point of the separation
coordinates}.
\end{Definition}

\begin{figure}[ht]
\centering
\begin{tikzpicture}[x=.85cm,y=.85cm,>=Latex,
  every node/.style={font=\small}]
\foreach \panel/\shift in {0/-3.9,1/3.9}{
  \begin{scope}[xshift=\shift cm]
    \fill[green!5] (0,0) circle (1.55);
    \begin{scope}
      \clip (0,0) circle (1.55);
      \fill[blue!9] (-1.6,-1.6) rectangle (0,1.6);
      \ifnum\panel=0
        \clip (-1.6,-1.6) rectangle (0,1.6);
      \fi
      \foreach \a in {-1.6,-1.2,-.8,-.4,0,.4,.8,1.2,1.6}
        \draw[gray!70] plot[domain=-1.8:1.8,samples=35]
          ({\a+.08*\x*\x},\x);
      \foreach \b in {-1.6,-1.2,-.8,-.4,0,.4,.8,1.2,1.6}
        \draw[gray!70] plot[domain=-1.8:1.8,samples=35]
          (\x,{\b+.06*\x*\x});
    \end{scope}
    \draw[green!45!black,thick] (0,0) circle (1.55);
    \draw[blue!60!black,dashed] (0,-1.55)--(0,1.55);
    \node[blue!60!black] at (-.9,1.1) {$\Omega\cap V$};
    \node[green!45!black] at (1.15,1.3) {$V$};
    \fill (0,0) circle (1.7pt);
    \node[anchor=north west,fill=white,inner sep=1pt]
      at (.05,-.05) {$P$};
    \ifnum\panel=0
      \node[align=center] at (0,-2.0)
        {coordinates $(x,y)$ on $\Omega\cap V$};
    \else
      \node[align=center] at (0,-2.0)
        {regular extension of $(\widetilde x,\widetilde y)$ to $V$};
    \fi
  \end{scope}
}
\draw[->,thick] (-1.7,0)--(1.7,0);
\node[align=center] at (0,.65)
  {$\widetilde x=\varphi(x)$\\$\widetilde y=\psi(y)$};
\end{tikzpicture}
\caption{}
\end{figure}

\begin{Proposition}\label{prop:singular-point-equivalence}
Let \((M^2, g)\) be a real-analytic
(pseudo-)Riemannian surface.  Let \((x,y)\) denote a system of
regular real or complex orthogonal separation coordinates on an open
set \(\Omega\subset M^2\), let \(P\in\overline{\Omega}\), and let
\(F\) be the quadratic integral associated with this separation.
Suppose that \(F\) extends real analytically to a neighbourhood \(V\)
of \(P\). Then \(P\) is a singular point of the separation coordinates if and
only if
\begin{equation*}
   D_L(P):=\left.
\bigl((\operatorname{tr}L)^2-4\det L\bigr)
\right|_P
=
0,
\end{equation*}
where \(L\) is the operator defined by~\eqref{operator} on \(V\).
\end{Proposition}

Thus, the singular points of the separation coordinates associated with
\(F\) are precisely the points at which the eigenvalues of \(L\)
coincide. The importance of this characterization is that it defines
the singular points invariantly.

Indeed, separation of variables gives rise to a quadratic first
integral \(F\) and hence to the associated operator \(L\). Its trace and determinant do not depend on the choice
of local coordinates. Let
\[
\Delta_L
:=
\left\{
P\in M^2:D_L(P)=0
\right\}
\]
be the singular set of the separation coordinates associated with
\(F\).

The set \(\Delta_L\) is invariant under changes of local coordinates.
It may, however, depend on the choice of the quadratic integral:
different quadratic integrals may determine inequivalent systems of
separation coordinates and hence different singular sets.
The replacement $\widetilde F=aF+bH$, $a\ne0$, preserves the
singular set, since $\widetilde L=aL+\frac b2\Id$ and
$D_{\widetilde L}=a^2D_L$.

Our aim is to describe all possible local singularities of the
separation coordinates, and, in particular, the possible local structures
of the set \(\Delta_L\).

\begin{Theorem}\label{thm:singular-set}
Let \(g\) be a
real-analytic (pseudo-)Riemannian metric on an open disk
\(D^2\). Suppose that, in a neighborhood \(U\) of a point
\(Q\in D^2\), its geodesic flow admits real or complex orthogonal separation of variables. Then
the following assertions hold 

\begin{itemize}
\item the geodesic flow is quadratically integrable on the whole disk,
and the associated operator is real analytic;

\item real or complex orthogonal separation of
variables is possible almost everywhere on the disk;

\item for every point \(P\in D^2\), there exists a neighborhood \(V\) such that \(V\cap\Delta_L\) is either empty, consists of a single point, is homeomorphic to a line, or is homeomorphic to the union of two intersecting lines.
\end{itemize}
\end{Theorem}

We are now ready to describe the local forms of systems of
separation coordinates near singular points.
Replacing $L$ by $L-\lambda_0\Id$, where $\lambda_0$ is
the repeated eigenvalue at the singular point $P$,
we assume that both eigenvalues of $L(P)$ are zero.
This replacement does not change the eigendirections.

\begin{Theorem}\label{thm:separation-coordinates}
Near every singular point the eigendirections of \(L\) determine a system of real
or complex orthogonal separation coordinates, up to admissible
transformations.

The following formulas describe the separation coordinates
$(\lambda_1,\lambda_2)$, up to admissible transformations

\begin{itemize}
    \item If $\det L \equiv0$
    
    \begin{itemize}
        \item $(\lambda_1, \lambda_2) = \left( xy, \frac{x}{y} \right)$;
        \item $(\lambda_1, \lambda_2) = \left(x^2 + y^2, \arctan(\frac{y}{x})\right)$;
        \item $(\lambda_1, \lambda_2) = \left(y, x - \ln|y|\right)$;
        \item $(\lambda_1, \lambda_2) = \left(y, x + \dfrac{1}{(2k-2)y^{2k-2}}\right)$, \ $k \geq 2$;
        \item $(\lambda_1,\lambda_2)=
\left(y,x+\frac{\varepsilon}{(2k-1)y^{2k-1}}\right),
\ k\ge1,\ \varepsilon=\pm1$
    \end{itemize}
    \item If $\det L \not \equiv 0$, then $\lambda_1$ and $\lambda_2$ are roots of 
    $$
    \lambda^2 - \tr L \lambda + \det L = 0,
    $$
   where \(\operatorname{tr}L\) and \(\det L\) are given in
    Table~\ref{tab:singularities}.
\end{itemize}
\end{Theorem}

\begin{table}[H]
\centering
\small
\caption{}
\label{tab:singularities}

\renewcommand{\arraystretch}{1.35}
\setlength{\tabcolsep}{4pt}
\begin{tabularx}{\textwidth}{
Y|
>{\centering\arraybackslash}p{2.10cm}|
Y}
\(\boldsymbol{\operatorname{tr}L}\) &
\(\boldsymbol{\det L}\) &
\textbf{Parameters} \\
\hline

\(2x\) & \(-y^2\) & \\ \hline

\(x\) & \(-y\) & \\ \hline

\(\displaystyle
\alpha xy^{2k-1}
+y^k\sum_{i=0}^{k-1}c_i y^i\)
&
\(-\varepsilon y^d\)
&
\(\begin{gathered}
k\geq1,\quad d\geq2k+1,\quad \varepsilon=\pm1,\\[-1mm]
\varepsilon=1\ \text{if }d\text{ is odd},\\[-1mm]
\alpha=kc_0^2\left(1-\dfrac{k}{d}\right)\neq0
\end{gathered}\) \\ \hline

\(\displaystyle
\alpha xy^s
+y^{s-k+1}\sum_{i=0}^{k-1}c_i y^i
+2\varepsilon y^k\)
&
\(y^{2k}\)
&
\(\begin{gathered}
k\geq1,\quad s\geq2k,\quad \varepsilon=\pm1,\\[-1mm]
\alpha=2\varepsilon kc_0\neq0
\end{gathered}\) \\ \hline

\(\displaystyle
\alpha xy^{2k-1}
+y^k\sum_{i=0}^{k-1}c_i y^i\)
&
\(-\varepsilon y^{2k}\)
&
\(\begin{gathered}
k\geq1,\quad \varepsilon=\pm1,\\[-1mm]
\alpha=\dfrac{k}{2}(c_0^2+4\varepsilon)\neq0
\end{gathered}\) \\ \hline

\(\displaystyle
\alpha xy^{2k}
+y^{k+1}\sum_{i=0}^{k-1}c_i y^i\)
&
\(-y^{2k+1}\)
&
\(k\geq1\), \ $\alpha = 2k+1$ \\ \hline

\(2x\) & \(y^2\) & \\ \hline

\(2xy^{s-1}-2y\) & \(y^2\) & \(s\geq2\)
\end{tabularx}
\end{table}

\section{Proofs}

\subsection{Proof of Theorem \ref{thm:Killing-separation}}

Real or complex nonorthogonal separation implies
the local existence of a nontrivial Killing vector field;
see~\cite[Theorems~4.2 and~5.1]{Benenti1991}.

Let $K$ be Killing vector and $J=K^ip_i$ is a linear first integral, $F=J^2$ is a quadratic first integral. The associated
operator $L$ has eigenvalues $0$ and $g(K,K)$.
Where $g(K,K)\ne0$, these eigenvalues are real and distinct,
so Theorem~\ref{thm:stackel} gives real
orthogonal separation. If $g(K,K)\not\equiv0$, its zero
set has measure zero by analyticity. This proves the first
assertion.

Suppose now that $g(K,K)\equiv0$. The metric is necessarily
Lorentzian. In local analytic null coordinates,
\[
g=2\lambda(u,v)\,\dd u\,\dd v,\qquad \lambda\ne0.
\]
By analyticity and the null condition, after interchanging
$u$ and $v$ if necessary, we have $K=f(u,v)\partial_u$.
The Killing equation gives
\[
\partial_vf=0,\qquad \partial_u(\lambda f)=0.
\]
Thus $f=f(u)$ and $\lambda(u,v)f(u)=a(v)$.
If $K$ vanished at a point with $u=u_0$, then $f(u_0)=0$
would imply $a\equiv0$, hence $f\equiv0$.
By analyticity and connectedness, this would contradict
$K\not\equiv0$. Therefore $K$ has no zeros.

Choose analytic coordinates $(s,t)$ with $K=\partial_t$.
The Killing and null conditions imply
\[
g=A(s)\,\dd s^2+2B(s)\,\dd s\,\dd t,\qquad B(s)\ne0.
\]
The analytic change of coordinates $x=\int^s B(\sigma)\,\dd\sigma,\
y=t+\int^s\frac{A(\sigma)}{2B(\sigma)}\,\dd\sigma$
has nonzero Jacobian and gives
\[
g=2\,\dd x\,\dd y,\qquad K=\partial_y.
\]
Consequently, $H=p_xp_y$, and the function $W(x,y;c_1,c_2)=\frac{c_1}{c_2}x+c_2y, \ c_2 \neq 0$
satisfies
\[
H(x,y, \frac{\partial W}{\partial x}, \frac{\partial W}{\partial y})=c_1,\qquad
\det\frac{\partial^2W}
{\partial(x,y)\,\partial(c_1,c_2)}
=\frac1{c_2}\ne0.
\]
The separation is nonorthogonal because
$g(\partial_x,\partial_y)=1$.

\subsection{Proof of Theorem \ref{thm:stackel}}

Suppose that the geodesic flow admits real or complex orthogonal separation near $P$.
By~\cite[Theorem~3.1]{Benenti1991} and, in the
complex-conjugate case,~\cite{degiovanni_rastelli},
there exists a real quadratic first integral $F$ such that,
in the separation coordinates,
\[
H=h_1p_1^2+h_2p_2^2,\qquad
F=f_1p_1^2+f_2p_2^2,\qquad
\delta:=h_1f_2-h_2f_1\ne0.
\]
Since $g$ is nondegenerate, $h_1h_2\ne0$.
By~\eqref{operator},
\[
L=
\begin{bmatrix}
\dfrac{f_2}{2h_2}&0\\[4pt]
0&\dfrac{f_1}{2h_1}
\end{bmatrix}.
\]
Its eigenvalues are distinct, since $\frac{f_2}{2h_2}-\frac{f_1}{2h_1}
=\frac{\delta}{2h_1h_2}\ne0.$

Suppose that $F$ is a real quadratic first integral and that
its associated operator $L$ has two distinct eigenvalues at $P$.
Since $L$ is Nijenhuis, there exist local real or
complex-conjugate coordinates $(q^1,q^2)$ in which
\[
L=\operatorname{diag}\bigl(\lambda_1(q^1),\lambda_2(q^2)\bigr),
\qquad \lambda_1-\lambda_2\ne0.
\]
Self-adjointness makes these coordinates orthogonal.
By~\eqref{operator}, we may therefore write
\[
H=h_1p_1^2+h_2p_2^2,\qquad
F=2\lambda_2h_1p_1^2+2\lambda_1h_2p_2^2.
\]
The coefficients of $p_1^2p_2$ and $p_1p_2^2$ in
$\{H,F\}=0$ give
\[
\partial_2\bigl((\lambda_1-\lambda_2)h_1\bigr)=0,
\qquad
\partial_1\bigl((\lambda_1-\lambda_2)h_2\bigr)=0.
\]
Hence
\[
H=\frac{a(q^1)p_1^2+b(q^2)p_2^2}{\lambda_1-\lambda_2},
\qquad
F=
2\frac{\lambda_2a(q^1)p_1^2+
       \lambda_1b(q^2)p_2^2}{\lambda_1-\lambda_2},
\]
where $ab\ne0$.

The function $W(q^1, q^2, c_1, c_2) = \int^{q^1}\sqrt{\frac{c_1\lambda_1(s)-c_2/2}{a(s)}}\,\dd s
+
\int^{q^2}\sqrt{\frac{c_2/2-c_1\lambda_2(t)}{b(t)}}\,\dd t$ is a complete integral. Indeed, direct calculations show that $H(q^1, q^2, \frac{\partial W}{\partial q^1}, \frac{\partial W}{\partial q^2}) = c_1$ and 

\[
\det\frac{\partial^2W}
{\partial(q^1,q^2)\,\partial(c_1,c_2)}
=\frac{\lambda_1-\lambda_2}{8abp_1p_2}\ne0.
\]

\subsection{Proof of Theorem \ref{main1}}

Let \(g\) be a quadratically integrable metric with quadratic integral
\(F\).  The required local normal forms of the pairs \((g,F)\) in a neighbourhood of a scalar singular point are
given in the Riemannian case by \cite{kol,bmf} and in the
pseudo-Riemannian case by \cite[Theorem~3]{A-B}.  We compute the
associated operator for each of these normal forms.

As we already noted  $L$ from \eqref{operator} is a Nijenhuis operator (see \cite{Akpan1, A-B}).
The proof of the Theorem proceeds as follows. We take a pair \((g,F)\) from the works cited above, construct the operator using~\eqref{operator}, and then perform an explicitly specified coordinate change \((x,y)\mapsto(u,v)\). We subsequently rewrite the invariants (the coefficients of the characteristic polynomial of \(L\)) in the new coordinates. Namely,
$$
\tr L = f(u,v), \quad \det L = -h(u,v).
$$

By~\cite[Corollaries~2.1--2.2]{BKM-1}, when the invariants
are functionally independent on a dense open set, they determine
$L$ uniquely through the formula
\begin{equation} \label{oper form}
L = \begin{bmatrix}
    f_u & f_v \\
    h_u & h_v
\end{bmatrix}^{-1} \begin{bmatrix}
    f & 1 \\
    h & 0
\end{bmatrix}  \begin{bmatrix}
    f_u & f_v \\
    h_u & h_v
\end{bmatrix}.
\end{equation}
This formula is used only where the invariant Jacobian is
invertible, and the resulting analytic expression extends across
its zero set. When the invariants are functionally dependent,
we transform the operator directly by
$L_{\mathrm{new}}=D\Phi\,L_{\mathrm{old}}\,(D\Phi)^{-1}$.

We begin with the pseudo-Riemannian case and use the notation of Theorem~3 in~\cite{A-B}. There, the metric is written as \(g=2\lambda(x,y)\dd x\dd y\). The condition that the quadratic integral \(F\) Poisson-commutes with the Hamiltonian \(H\) implies that \(F=a(x)p_x^2+b(x,y)p_xp_y+c(y)p_y^2\). Hence, in these coordinates, the operator \(L\) defined by~\eqref{operator} has the following form:
\begin{equation} \label{oper matrix}
    L = \begin{bmatrix}
        \frac{\lambda b}{2} & -a\lambda \\
        -c \lambda & \frac{\lambda b}{2}
    \end{bmatrix}.
\end{equation}

Without loss of generality, we assume that $P=(0,0)$ and
$L(P)=0\cdot\operatorname{Id}$; otherwise, we replace $L$ by
$L-\lambda_0\operatorname{Id}$.

\textbf{Pseudo-Riemannian Case}

\textbf{C1.1}

In this case, the metric and quadratic integral have the form
 \begin{equation*}
        \begin{aligned}                           
        g &= \frac{X\left(\sqrt{x} + \sqrt{y}\right) - X\left(\sqrt{x} - \sqrt{y}\right)}{\sqrt{x}\sqrt{y}}\, \dd x \dd y, \\
         F &= x\, p^2_x  + y\, p^2_y - 2\sqrt{x}\sqrt{y} \, \frac{X\left(\sqrt{x} + \sqrt{y}\right) + X\left(\sqrt{x} - \sqrt{y}\right)}{X\left(\sqrt{x} + \sqrt{y}\right) - X\left(\sqrt{x} - \sqrt{y}\right)}\,p_xp_y, 
        \end{aligned}
        \end{equation*} 
where $X(t)$ is even and analytic near zero and satisfies
$X''(0)\neq0$. Without loss of generality, we may assume that
$X(0)=0$. Otherwise, we replace $X(t)$ by $X(t)-X(0)$; the metric
remains unchanged, while the quadratic integral becomes $ F \longrightarrow F + \operatorname{const} \cdot H,$
where $H$ is the Hamiltonian of the geodesic flow. This corresponds
exactly to shifting $L$ so that $L(P)=0\cdot\operatorname{Id}$.

The operator \(L\) defined by~\eqref{oper matrix} has the form
$$
L = \begin{bmatrix}
    - \frac{X(\sqrt{x} + \sqrt{y}) + X(\sqrt{x} - \sqrt{y})}{2} & -\sqrt{\frac{x}{y}} \frac{X(\sqrt{x} + \sqrt{y}) - X(\sqrt{x} - \sqrt{y})}{2} \\
     -\sqrt{\frac{y}{x}} \frac{X(\sqrt{x} + \sqrt{y}) - X(\sqrt{x} - \sqrt{y})}{2} & - \frac{X(\sqrt{x} + \sqrt{y}) + X(\sqrt{x} - \sqrt{y})}{2}
\end{bmatrix}.
$$

Denote \(X_{\pm}=X(\sqrt{x}\pm\sqrt{y})\) and set $u=-\frac{X_++X_-}{2},\ v=\sqrt{X_+X_-}.$ Here and below the square root is the signed analytic branch, not the nonnegative principal square root. To see that this branch exists, use \(X(0)=0\), \(X''(0)\ne0\), and the evenness of \(X\). Hadamard's lemma gives \(X(t)=t^2Z(t^2)\), where \(Z\) is analytic and \(Z(0)\ne0\). Since $Z_+Z_-$ is analytic and has value $Z^2(0)>0$ at the
origin, it has a real-analytic square root.

$$
u = - \frac{X_+ + X_-}{2} = -\sum_{i = 1}^{\infty} X_{2i} \sum_{k = 0}^i \binom{2i}{2k} x^{i-k}y^k, \quad v = (x-y)\sqrt{Z_+ Z_-},
$$
where $X_j$ are the Taylor coefficients of $X(t)$ at zero.
The Jacobian at the origin is $2Z(0)\sqrt{Z(0)^2}\ne0$.
In coordinates $(u,v)$, the trace and determinant of $L$ are
$$
\tr L = 2u, \quad \det L = v^2.
$$

Formula~\eqref{oper form} then gives case (SS1) $L = \begin{bmatrix}
    2u & -v \\ 
    v & 0
\end{bmatrix}.$

\textbf{C1.2, C1.3, C2}

We first consider case C1.3. For case C1.3, the metric and quadratic integral have the form
 \begin{equation*}
        \begin{aligned}            
        g &= \frac{
        X\left(y\bigl(1 +  2\sqrt{x}y^{\frac{s-2}{2}}\bigr)^{\frac{2}{2-s}}\right) - 
        X\left(y\bigl(1 -  2\sqrt{x}y^{\frac{s-2}{2}}\bigr)^{\frac{2}{2-s}}\right)}
        { \frac{2}{2-s}\sqrt{x} \, y^{\frac{s}{2}}} \, \dd x \dd y, \\ 
        F &= xp^2_x 
        + \left(\tfrac{2}{2-s}\right)^2y^sp^2_y
        - \tfrac{4\sqrt{x}}{2-s} \, y^{\frac{s}{2}} \frac{
        X\left(y \bigl( 1 {+}  2\sqrt{x}y^{\frac{s-2}{2}}\bigr)^{\frac{2}{2-s}}\right) + 
        X\left( y\bigl(1 {-}  2\sqrt{x}y^{\frac{s-2}{2}}\bigr)^{\frac{2}{2-s}}\right)}
        {X\left(y\bigl(1 {+}  2\sqrt{x}y^{\frac{s-2}{2}}\bigr)^{\frac{2}{2-s}}\right) - 
        X\left(y\bigl(1 {-}  2\sqrt{x}y^{\frac{s-2}{2}}\bigr)^{\frac{2}{2-s}}\right)}p_xp_y , 
        \end{aligned}
        \end{equation*} 
where $X(t)$ is analytic near zero, $X'(0)\neq0$, and $s>2$.
As in the previous case, we may assume that $X(0)=0$; otherwise, we
replace $X(t)$ by $X(t)-X(0)$. 

In this case, the operator \(L\) defined by~\eqref{oper matrix} has the form
$$
L = \begin{bmatrix}
    -\frac{X_+ + X_-}{2} & -\frac{\sqrt{x}}{2} \frac{2-s}{2} y^{-\frac{s}{2}} (X_+ - X_-) \\
    -\frac{y^{\frac{s}{2}}}{\sqrt{x}(2-s)} (X_+ - X_-) & -\frac{X_+ + X_-}{2}
\end{bmatrix},
$$
where 
$X_{\pm} = X(y(1 \pm 2\sqrt{x} y^{\frac{s-2}{2}})^{\frac{2}{2-s}})$, $s > 2$.

Set $q=2/(2-s)$, where $s>2$ is an integer, and write
$X(t)=tY(t)$, with $Y(0)=X'(0)\ne0$.
Temporarily put $\beta=2\sqrt{x}\,y^{(s-2)/2}$ and define
\[
a_\pm=(1\pm\beta)^qY\bigl(y(1\pm\beta)^q\bigr),
\qquad X_\pm=ya_\pm.
\]
Choose the analytic branch $R$ satisfying
$R^2=a_+a_-$ and $R(0,0)=Y(0)$, and set
\[
v=yR,\qquad
u=\frac{2v-X_+-X_-}{2v^{s-1}}.
\]
The product $a_+a_-$ and the sum $a_++a_-$ are even in $\beta$.
Moreover,
\[
a_+-a_-=\beta T,\qquad T(0,0)=2qY(0)\ne0,
\]
where $T$ is analytic and even in $\beta$.
Thus $R$ and $T$ are analytic functions of $(y,\beta^2)$,
and hence of $(x,y)$, since $\beta^2=4xy^{s-2}$.
Using
\[
2v-X_+-X_-
=-y\,\frac{(a_+-a_-)^2}{a_++a_-+2R},
\]
we obtain the analytic expression $u=-\frac{2xT^2}{R^{s-1}(a_++a_-+2R)}.$
All denominators are nonzero at the origin. In particular,
\[
v_y(0,0)=Y(0),\qquad
u_x(0,0)=-2q^2Y(0)^{2-s}\ne0,\qquad v_x(0,0)=0.
\]
Therefore $(u,v)$ are local analytic coordinates. In these coordinates, the trace and determinant are
$$
\tr L = 2 u v^{s-1} - 2v, \quad \det L = v^2,
$$
and formula~\eqref{oper form} gives case (SS2)
$$
L = \begin{bmatrix}
    -v(1 + (s-3)uv^{s-2}) & -u(2 + (s-1)(s-3)uv^{s-2}) \\
    v^s & -v(1 - (s-1)uv^{s-2})
\end{bmatrix}.
$$

The remaining cases C1.2 and C2 can be treated together. Put
$s=2k$ and let $\rho(t,y)$ be the analytic flow of
\[
f_0(y)\partial_y,
\qquad
f_0(y)=
\begin{cases}
y,&k=1\quad\text{(C1.2)},\\[2pt]
\dfrac{y^k}{1+y^{k-1}},&k\ge2\quad\text{(C2)}.
\end{cases}
\]
Thus $\rho(0,y)=y$; in C1.2, $\rho(t,y)=ye^t$.
After subtracting the constant energy term, assume $X(0)=0$
and write $X(t)=tY(t)$, with $Y(0)\ne0$.
The normal forms in~\cite[Theorem~3]{A-B} give
\[
\tr L=-X_+-X_-,\qquad \det L=X_+X_-,
\qquad X_\pm=X\bigl(\rho(\pm2\sqrt{x},y)\bigr).
\]
To justify the coordinate change analytically, first regard $t$ and
$y$ as independent variables and set
\[
a_\pm(t,y)=\frac{X(\rho(\pm t,y))}{y}.
\]
These functions are analytic, since $\rho(t,0)=0$.
Moreover,
\[
a_+-a_-=t y^{k-1}T(t^2,y),\qquad T(0,0)=2Y(0),
\]
with $T$ analytic. Indeed, $\partial_tX(\rho(t,y))
=X'(\rho(t,y))f_0(\rho(t,y))$ is divisible by $y^k$;
the difference $X(\rho(t,y))-X(\rho(-t,y))$ is odd in $t$
and its leading term is $2t y^kY(0)$.
The product $a_+a_-$ and the sum $a_++a_-$ are even in $t$.
Choose the even analytic square root $R$ with
$R^2=a_+a_-$ and $R(0,0)=Y(0)$.
After substituting $t^2=4x$, define
\[
v=yR,\qquad
u=\frac{2v-X_+-X_-}{2v^{2k-1}}
=-\frac{2xT^2}{R^{2k-1}(a_++a_-+2R)}.
\]
The last expression is analytic and has nonzero denominators near
the origin. Furthermore,
\[
v_y(0,0)=Y(0),\qquad v_x(0,0)=0,\qquad
u_x(0,0)=-2Y(0)^{2-2k}\ne0.
\]
Thus $(u,v)$ are analytic coordinates and
\[
\tr L=2uv^{2k-1}-2v,\qquad \det L=v^2.
\]
Formula~\eqref{oper form} gives SS2 with $s=2k$, including $s=2$.

\textbf{C3}
After subtracting the constant energy term, the pair is
\[
g=h(xy)\,\dd x\,\dd y,\qquad
F=(xp_x-yp_y)^2,\qquad h(0)\ne0.
\]
Formula~\eqref{oper matrix} gives $L=-\frac{h(xy)}2
\begin{bmatrix}xy&x^2\\y^2&xy\end{bmatrix}.$
Put $\varepsilon=\operatorname{sign}h(0)$ and define $u=x\sqrt{\frac{|h(xy)|}{2}},\
v=-\varepsilon y\sqrt{\frac{|h(xy)|}{2}}.$ This is an analytic coordinate change with nonzero Jacobian at
the origin. Direct multiplication gives
\[
D\Phi\,L=
\begin{bmatrix}uv&u^2\\v^2&uv\end{bmatrix}D\Phi,
\qquad \Phi=(u,v).
\]
Hence the transformed operator is SS3.

\textbf{C4}
Here
\[
g=2\,\dd x\,\dd y,\qquad
F=yp_y^2-(x+C)p_xp_y,
\]
and therefore
\[
L=\begin{bmatrix}-(x+C)/2&0\\-y&-(x+C)/2\end{bmatrix}.
\]
Removing its constant scalar term $-(C/2)\Id$ and applying
$u=-x$, $v=y$ gives directly
\[
L_{\mathrm{new}}=
\begin{bmatrix}u/2&0\\v&u/2\end{bmatrix},
\]
which is SS4.

Thus, we have covered all pseudo-Riemannian cases. We now consider the Riemannian metrics.

\textbf{Riemannian Case}

In the Riemannian case, there are three types of normal forms
$(g,F)$; see \cite{bmf}. One corresponds to the
$\operatorname{gl}$-regular case and will not be discussed here. We
consider the remaining two cases.

\textbf{The first case is:}
$$
g = \frac{h(x + r) - h(x - r)}{2r} (\dd x^2 + \dd y^2), \quad F = - r\frac{h(x + r) + h(x - r)}{h(x + r) - h(x - r)}(p^2_x + p^2_y) + (xp^2_x + 2yp_xp_y - xp^2_y),
$$
where $r = \sqrt{x^2+y^2}$, $h(t)$ is real analytic in a neighborhood of zero, and $h'(0) > 0$. Without loss of generality, assume $h(0) = 0$.

Once again, we construct the operator \(L\) using~\eqref{operator}.
$$
L = \begin{bmatrix}
    -\frac{h_+ - h_-}{2r}\left(x + r \frac{h_+ + h_-}{h_+ - h_-} \right) & -y\frac{h_+ - h_-}{2r} \\
    -y\frac{h_+ - h_-}{2r} &  \frac{h_+ - h_-}{2r}\left(x - r \frac{h_+ + h_-}{h_+ - h_-} \right)
\end{bmatrix},
$$
where $h_{\pm}=h(x\pm r)$. The new coordinates are
$v=\sqrt{-h(x+r)h(x-r)}$ and
$u=-\frac12(h(x+r)+h(x-r))$. To see that $u$ and $v$ are real analytic,
write $h(t)=tH(t)$, where $H(0)>0$. Then
$v=y\sqrt{H(x-r)H(x+r)}$ is real analytic near the origin. Moreover, $u$ is
even in $r$ and can therefore be written as
$u=R(x,r^2)=R(x,x^2+y^2)$ for a real-analytic function $R$.
The Jacobian at the origin is $-H(0)^2\ne0$.

In the new coordinates $\tr L = 2 u, \ \det L = -v^2,$
and formula~\eqref{oper form} gives case (SS5) $L = \begin{bmatrix}
    2 u & v \\
    v & 0
\end{bmatrix}.$

\textbf{The second case is:}
$$
g = f(x^2+y^2) (\dd x^2 + \dd y^2), \quad F = (yp_x - xp_y)^2,
$$
where $f(t)$ is a positive real-analytic function.

The operator \(L\) defined by~\eqref{operator} has the form $L = \begin{bmatrix}
    x^2 f & xy f \\
    xyf & y^2f
\end{bmatrix}.$ For $\Phi=(u,v)=(x\sqrt f,y\sqrt f)$, a direct calculation gives $D\Phi\,L=
\begin{bmatrix}u^2&uv\\uv&v^2\end{bmatrix}D\Phi.$ Since $\det D\Phi(0,0)=f(0)>0$, this is an analytic coordinate
change and the transformed operator is SS6.

\subsection{Proof of Proposition \ref{prop:singular-point-equivalence}}

Suppose that $P$ is regular. In the extended orthogonal
separation coordinates, both $g$ and $F$ are diagonal.
The equation $\{H,F\}=0$ implies that the difference of
the eigenvalues of $L$ satisfies homogeneous linear
first-order equations with regular coefficients.
If this difference vanished at $P$, uniqueness would
force it to vanish throughout a connected neighbourhood,
contrary to the distinctness of the eigenvalues on the
original separation domain.

Conversely, if $L(P)$ has distinct eigenvalues, the
construction in the proof of Theorem~\ref{thm:stackel}
gives regular separation coordinates near $P$ adapted
to its eigendirections. On the original domain, they
have the same coordinate curves as the given coordinates,
so they differ only by an admissible transformation.
Thus $P$ is regular.

The assertion follows because the eigenvalues coincide
exactly when $D_L(P)=0$.

\subsection{Proof of Theorem \ref{thm:singular-set}}
We extend the local quadratic integral in five steps.

\begin{enumerate}
\item
Choose real-analytic coordinates on the disk in which
\(D^2=\{x^2+y^2<1\}\) and \(Q=(0,0)\).

\item
The equation $\{H,F\}=0$ for a quadratic integral
is equivalent to the Killing tensor equation.
By \cite[Theorem~8]{KruglikovMatveev2016}, with
$n=2$ and $d=2$, its prolongation defines a linear
connection $\mathfrak D$ on a rank-six bundle $\mathcal E$,
whose parallel sections correspond bijectively to
quadratic first integrals. An alternative construction is given in
\cite[Corollary~17]{GoverLeistner2019}, applied with
$r=2$ to the projective structure $[\nabla^g]$:
the parallel metric volume density identifies their
weighted Killing tensors with ordinary Killing tensors.

In local fibre coordinates, parallel sections $s = s(x,y)$ satisfy
\[
 \mathfrak D_xs:=\partial_xs-A_xs=0,\qquad
 \mathfrak D_ys:=\partial_ys-A_ys=0,
\]
with analytic matrices $A_x,A_y$.
Let $s_0$ be the local parallel section corresponding
to $F$ on $U$.

\item
For each $X=(x,y)$ in the disk,  parallel transport $s_0(0)$
along the radius $\gamma_X(t)=tX$.
In local fibre coordinates, the transport equation is
\[
 \frac{d}{dt}w(t,X)
 =
 \bigl(xA_x(tX)+yA_y(tX)\bigr)w(t,X),
 \qquad w(0,X)=s_0(0).
\]
The existence and uniqueness theorem for linear ODEs
gives transport along the entire compact segment
$0\leq t\leq1$. Define $\widetilde s(X)=w(1,X)$.

\item
We use analytic dependence of ODE solutions on parameters:
if $f(t,w,\lambda)$ and the initial value
$w_0(\lambda)$ are analytic, then the solution of
\[
 \dot w=f(t,w,\lambda),\qquad
 w(0,\lambda)=w_0(\lambda)
\]
depends analytically on $\lambda$ locally along
any interval on which it exists (see \cite[Theorem~4.2]{Teschl2012}).

Here the endpoint $X=(x,y)$ is the parameter.
The radii $\gamma_X(t)=tX$ vary analytically with $X$,
and hence
\[
 f(t,w,X)=\bigl(xA_x(tx,ty)+yA_y(tx,ty)\bigr)w
\]
is analytic in $(t,w,X)$.
The initial value $w_0(X)=s_0(0)$ is constant.
Thus $\widetilde s(X)=w(1,X)$ is analytic near every
endpoint $X$. 

\item
Choose a smaller disk $B\subset U$ centred at the origin.
On each radius in $B$, the original section $s_0$
satisfies the same transport equation with the same
initial value. Uniqueness gives $\widetilde s=s_0$
on $B$, so $\mathfrak D\widetilde s=0$ there.
Since $\mathfrak D$ and $\widetilde s$ are analytic,
$\mathfrak D\widetilde s$ is an analytic section.
The identity theorem on the connected disk implies
$\mathfrak D\widetilde s=0$ everywhere.

By the correspondence in
\cite[Theorem~8]{KruglikovMatveev2016},
$\widetilde s$ determines a global analytic quadratic
first integral $\widetilde F$.
Since $\widetilde s=s_0$ on $B$, we have
$\widetilde F=F$ there, and hence on the connected
original domain $U$ by analyticity.
\end{enumerate}

Finally, $dH\wedge d\widetilde F$ is analytic and
nonzero on the original regular domain.
Its nonvanishing set is therefore open and dense
in phase space. This proves the first assertion.

Thus \(F\), and consequently $L$ by \eqref{operator} are globally defined and real analytic on \(D^2\). The operator \(L\) has two distinct eigenvalues near \(Q\), because
\(F\) is associated with an orthogonal separation.
Hence $D_L\not\equiv0$
on \(D^2\). By Proposition~\ref{prop:singular-point-equivalence}, \(\Delta_L\) is
precisely the singular set of the orthogonal separation coordinates
associated with \(F\).

For every point \(P\in D^2\setminus\Delta_L\), the operator \(L\) has
two distinct eigenvalues in a sufficiently small neighborhood of
\(P\). Theorem~\ref{thm:stackel} therefore gives real or complex orthogonal separation there. It remains to describe the local structure of $\Delta_L$.

There are two cases at a point with repeated eigenvalues.
If $L(P)$ is non-scalar, it is a single Jordan block and is
$\operatorname{gl}$-regular in a neighbourhood of $P$. After subtracting
its constant eigenvalue at $P$ and applying a constant change of
basis, Theorem~5.1 of~\cite{BKM-3} applies. In its $L_{\mathrm{nd}}$
normal form the discriminant is $x^2+4y$, whose zero set is a smooth
curve. For the series $M,O,P,S$, its local zero set is $y=0$;
see also the description of the singular set in Section~5 of that
paper. The forms $L_{\mathrm{nil}}$ and $N$ have identically zero
discriminant and are excluded by $D_L\not\equiv0$ and analyticity.
Thus the non-scalar repeated-eigenvalue case gives a local line.

If $L(P)$ is scalar, $D_L\not\equiv0$ implies that $P$ is a singular
scalar point in the sense of Theorem~\ref{main1}. That theorem applies
because \(L\) is real analytic and is associated, by construction,
with the metric \(g\) and the global quadratic integral \(F\).
The discriminants of its six normal forms are
\[
\begin{array}{c|c|c|c|c|c}
\mathrm{SS1} & \mathrm{SS2} & \mathrm{SS3}
& \mathrm{SS4} & \mathrm{SS5} & \mathrm{SS6}\\
\hline
4(x^2-y^2)
& 4xy^s(xy^{s-2}-2)
& 4x^2y^2
& 0
& 4(x^2+y^2)
& (x^2+y^2)^2
\end{array}
\qquad (s\geq2).
\]
Adding the constant scalar term $\lambda_0\operatorname{Id}$ does
not change the discriminant. In SS2, the factor $xy^{s-2}-2$ is
nonzero near the origin. Hence SS1--SS3 give two intersecting lines,
whereas SS5 and SS6 give only the origin.

In case SS4, the discriminant vanishes identically in a neighborhood
of the point. This case cannot occur. Indeed, since \(D_L\) is real
analytic, its vanishing on a nonempty open subset of the connected
disk would imply $D_L\equiv0$
on \(D^2\), contradicting \(D_L\neq0\) near \(Q\).

Thus, near every point, \(\Delta_L\) has one of the stated local
topological types. Since it is the zero set of the nonzero real-analytic
function $D_L$, it has measure zero. Therefore,
real or complex orthogonal separation is possible
in a neighborhood of almost every point of \(D^2\).

\subsection{Proof of Theorem \ref{thm:separation-coordinates}}

By Proposition~\ref{prop:singular-point-equivalence}, the complement
of the singular set is precisely the set on which \(L\) has two
distinct eigenvalues.  By Theorem~\ref{thm:stackel}, the corresponding
eigendirections determine a regular system of real or complex
orthogonal separation coordinates, uniquely up to admissible
transformations.

For the normal forms with \(\det L\equiv0\), the coordinate pairs
listed above are obtained by direct integration of the
eigendirections.  More precisely, these are the cases SS3, SS6, and
the \(M\)-series (see ~\cite{BKM-3}), respectively.

\begin{figure}[H]   
    \centering
    \begin{minipage}{0.3\textwidth}
        \centering
        \includegraphics[width=5cm]{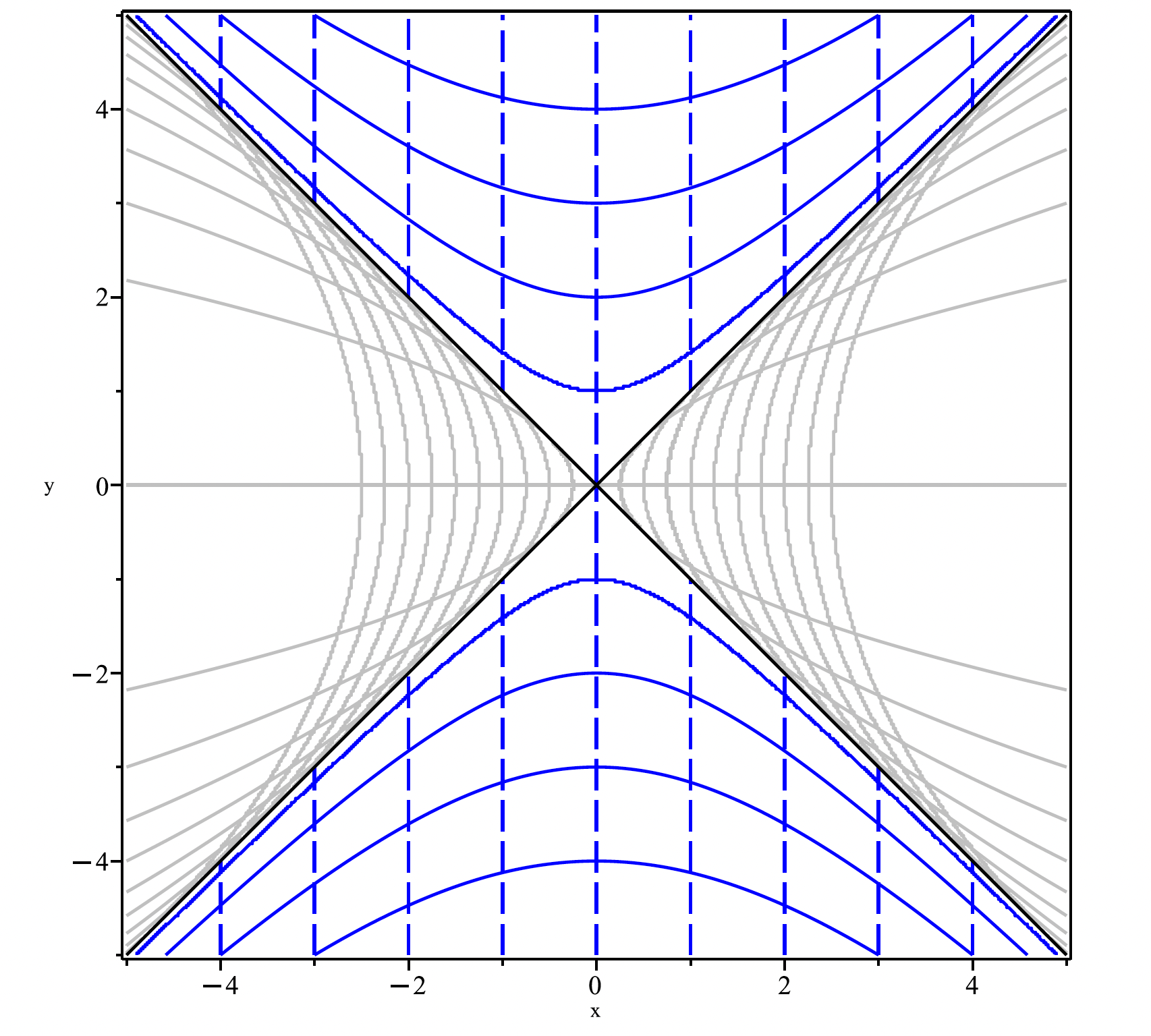}
        \caption{SS1}
        \label{fig:first}
    \end{minipage}
    \hfill
     \begin{minipage}{0.3\textwidth}
        \centering
        \includegraphics[width=5cm]{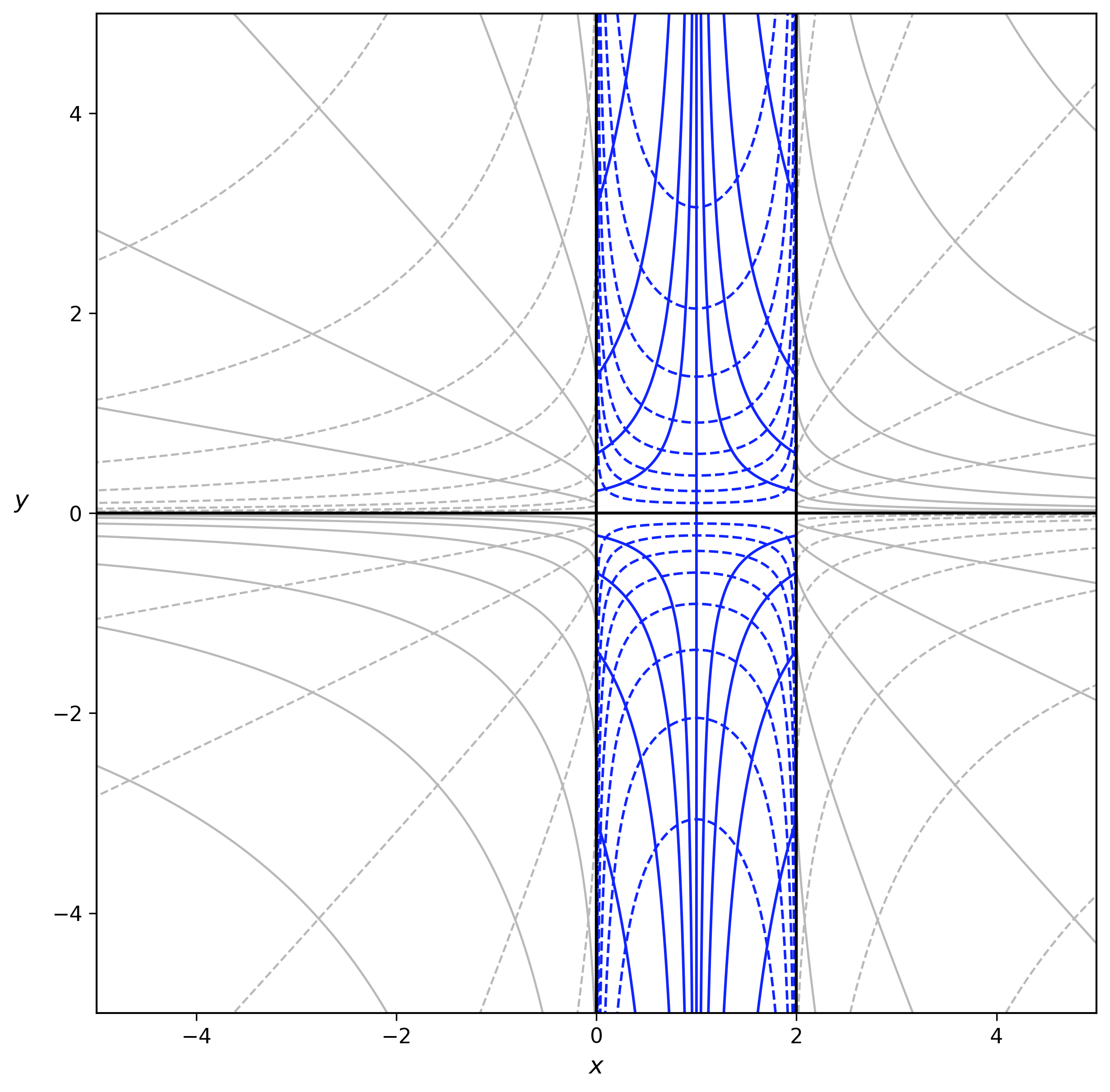}
        \caption{SS2, $s =2$}
    \end{minipage}
    \hfill
     \begin{minipage}{0.3\textwidth}
        \centering
        \includegraphics[width=5cm]{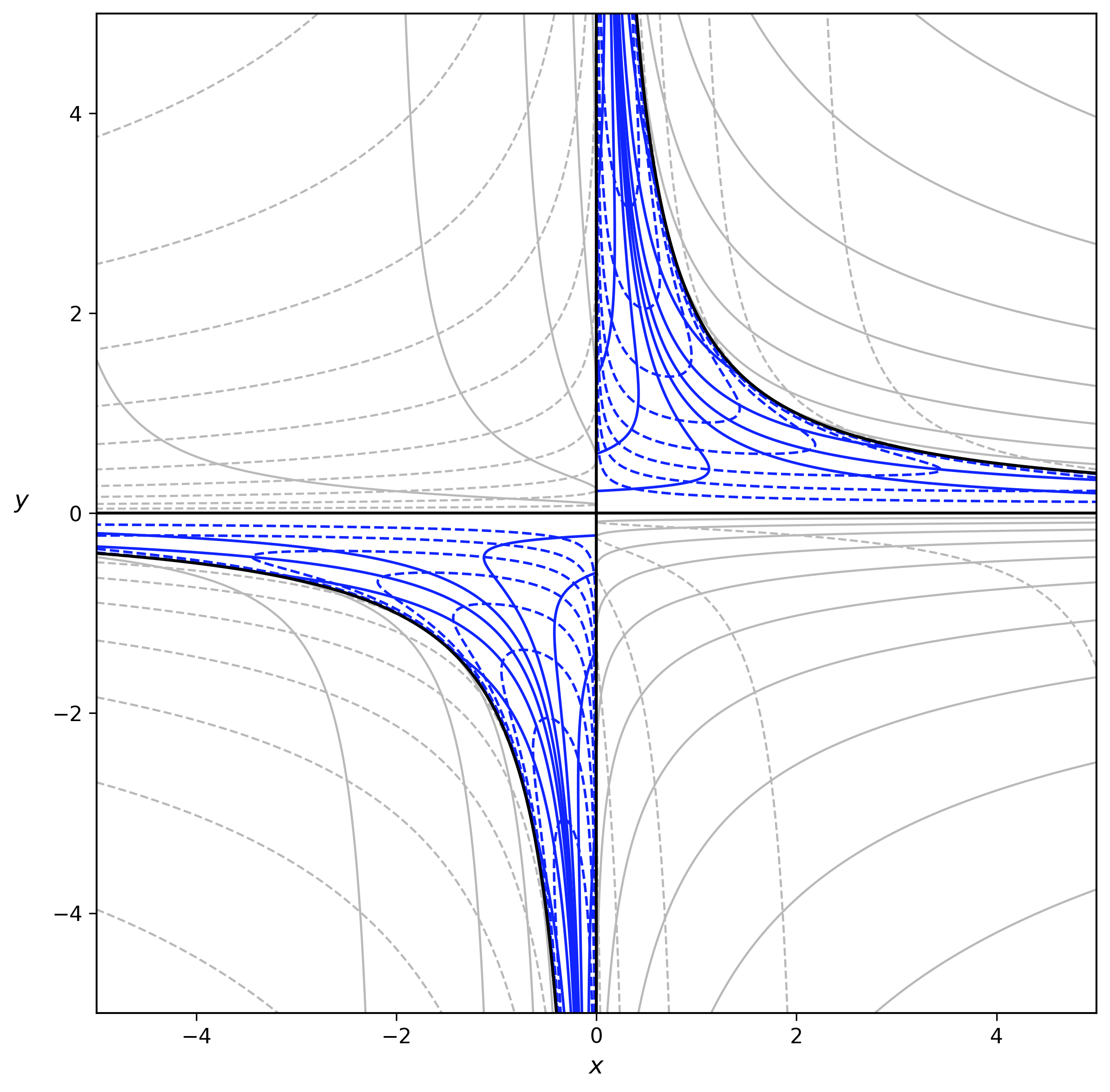}
        \caption{SS2, $s = 3$}
    \end{minipage}
\end{figure}

\begin{figure}[H]
  \begin{minipage}{0.3\textwidth}
        \centering
        \includegraphics[width=5cm]{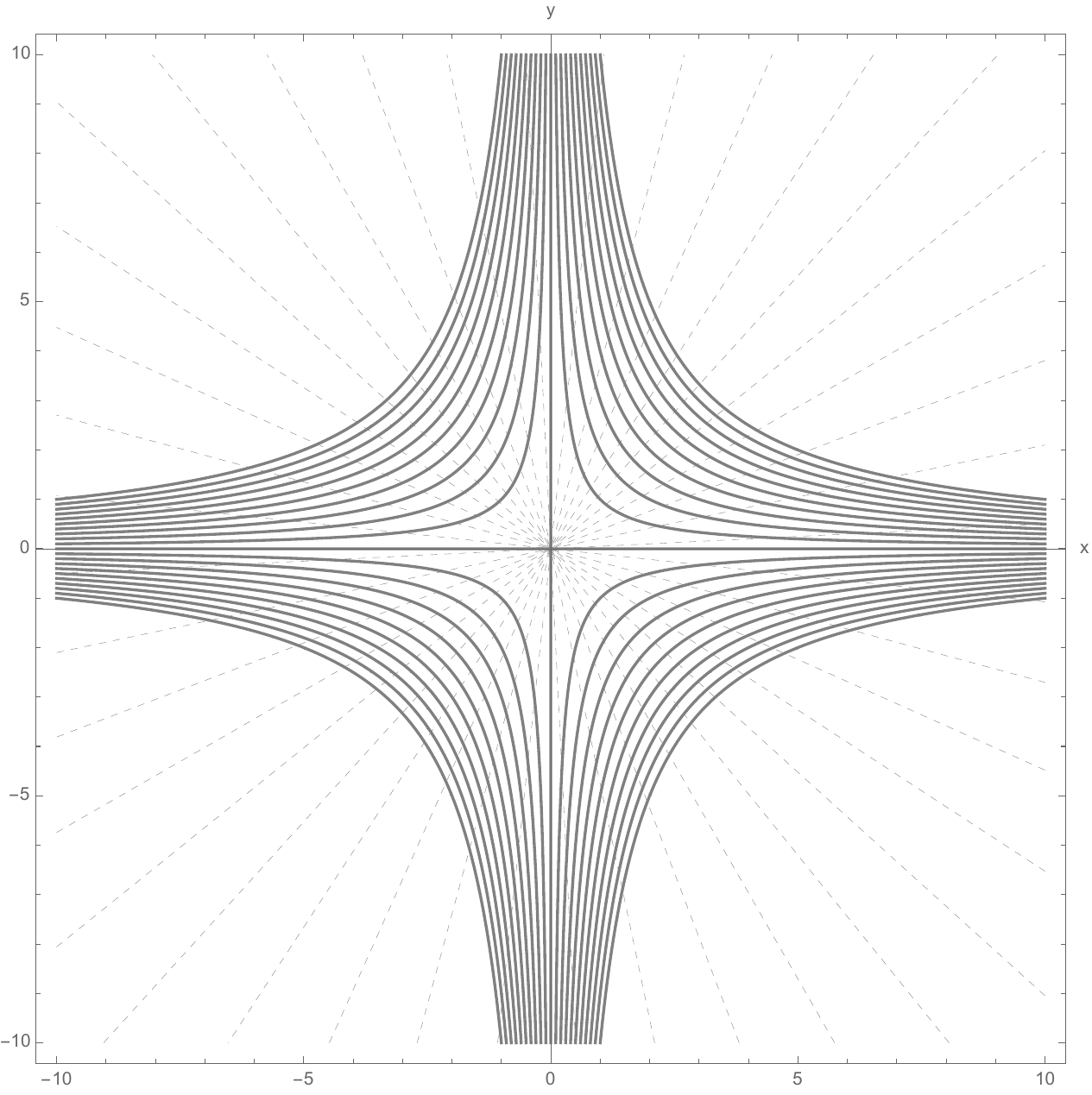}
        \caption{SS3}
        \label{fig:rindler}
    \end{minipage}
    \hfill
      \begin{minipage}{0.3\textwidth}
        \centering
        \includegraphics[width=5cm]{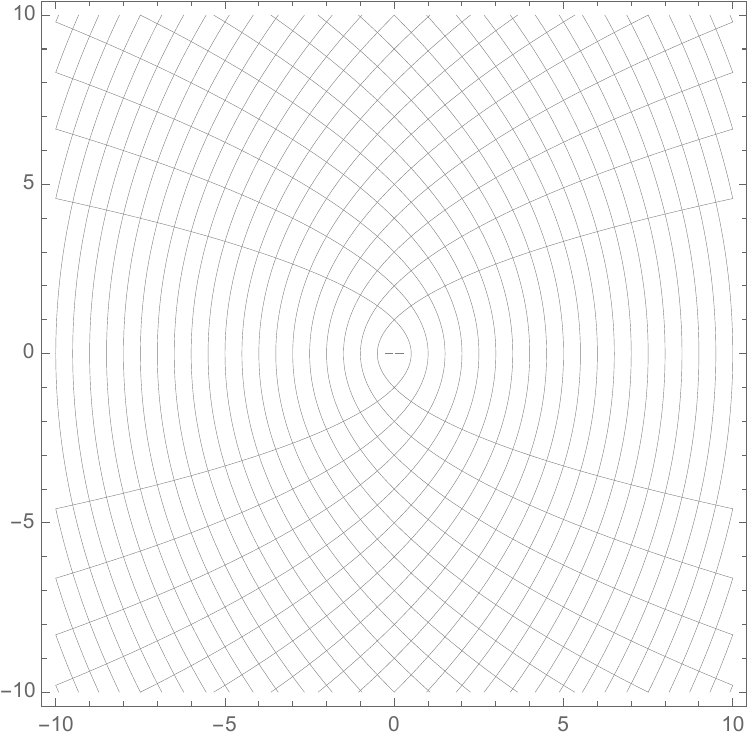}
        \caption{SS5}
        \label{fig:third}
    \end{minipage}
     \hfill
    \begin{minipage}{0.3\textwidth}
        \centering
        \includegraphics[width=5cm]{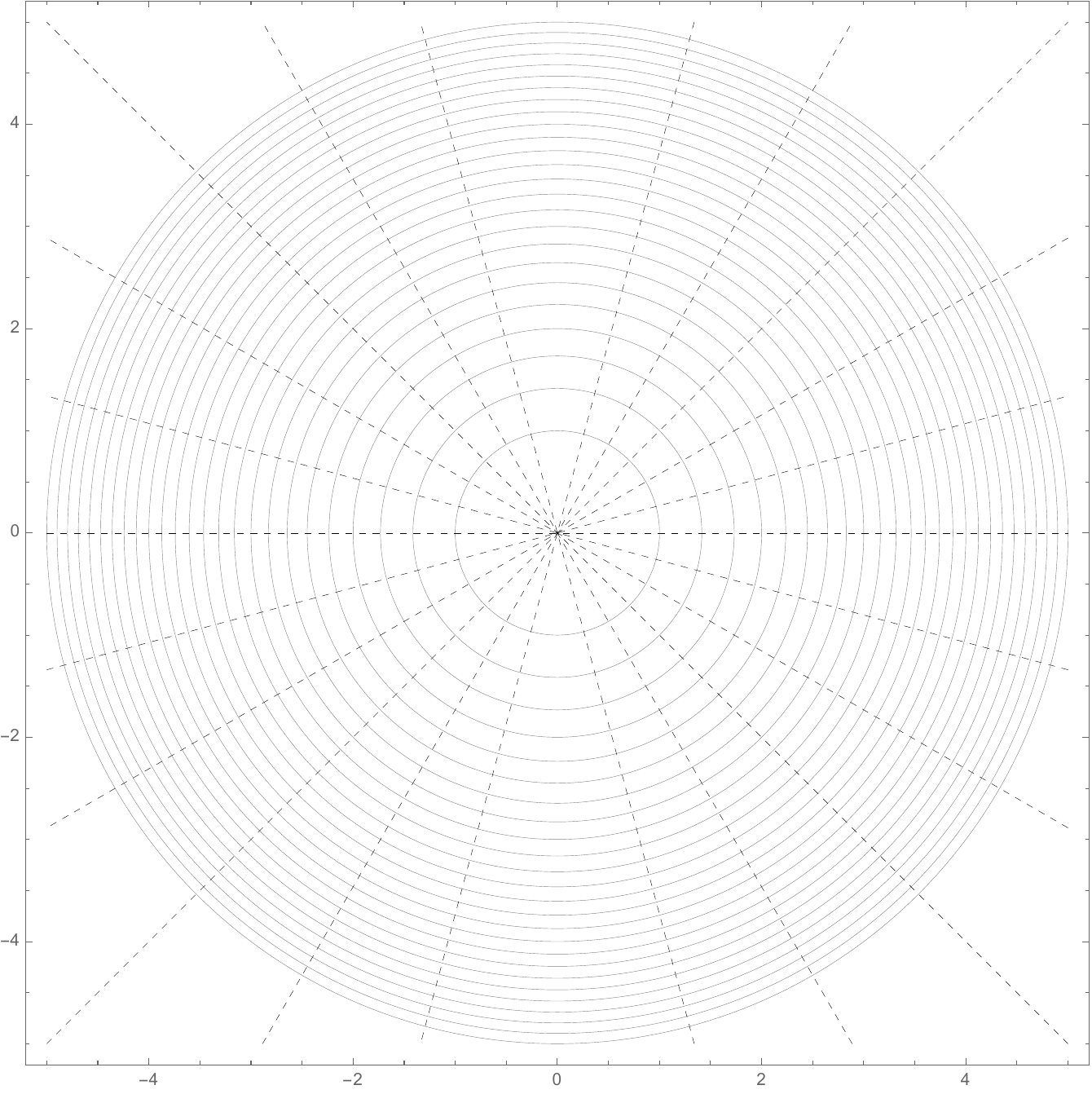}
        \caption{SS6}
        \label{fig:fourth}
    \end{minipage}
\end{figure}

\begin{Remark}
Figures~\ref{fig:first}--\ref{fig:fourth} illustrate the separation
coordinates from Theorem~\ref{thm:separation-coordinates}. Real
separation coordinates are shown in gray. Where separation is given
by nonreal complex-conjugate eigenvalues, the level curves of their
real and imaginary parts are shown in blue. 
\end{Remark}

\end{document}